\documentclass{article}
\usepackage{spconf,graphicx}
\usepackage[cmex10]{amsmath}
\usepackage{cite}
\usepackage{array}
\usepackage{mdwmath}
\usepackage{mdwtab}
\usepackage[caption=false,font=footnotesize]{subfig}
\usepackage{stfloats}
\usepackage{url}

\usepackage[american,fulldiode]{circuitikz} 
\usepackage{graphicx} 
\usepackage{multirow}
\usepackage{amsfonts}
\usepackage{bbm}
\usepackage{algorithm}
\usepackage{algpseudocode}
\usepackage{multirow}
\usepackage[section]{placeins}
\usepackage{epstopdf}

\usepackage{slashbox}
\usepackage{amssymb}

\usepackage{lipsum,amsmath,multicol}

\usepackage{multirow}

\usepackage{float}
\usepackage{tablefootnote}
\usepackage[multiple]{footmisc}

\title{Moment Relaxations of Optimal Power Flow Problems: \\ Beyond the Convex Hull}
\threeauthors
  {Daniel K. Molzahn\thanks{The authors gratefully acknowledge a discussion with \mbox{Dr.~Carleton} \mbox{Coffrin} as one of the impetuses for this paper. \newline\hspace*{15pt} Argonne National Laboratory's work was supported by the U.S. Department of Energy, Office of Electricity Delivery and Energy Reliability under contract DE-AC02-06CH11357.}}
	{Argonne National Laboratory\\
	Energy Systems Division\\
	dmolzahn@anl.gov}
  {C\'edric Josz}
	{R\'eseau  de  Transport d'Electricit\'e \\ R\&D Department\\ cedric.josz@rte-france.com}
  {Ian A. Hiskens}
  {University of Michigan \\
  EECS Department\\
  hiskens@umich.edu}
\begin{document}
%
\maketitle
\begin{abstract}
Optimal power flow (OPF) is one of the key electric power system optimization problems. ``Moment'' relaxations from the Lasserre hierarchy for polynomial optimization globally solve many OPF problems. Previous work illustrates the ability of higher-order moment relaxations to approach the convex hulls of OPF problems' non-convex feasible spaces. Using a small test case, this paper focuses on the ability of the moment relaxations to globally solve problems with objective functions that have unconstrained minima at infeasible points inside the convex hull of the non-convex constraints.
\end{abstract}
\begin{keywords}
Optimal power flow, convex relaxation
\end{keywords}
\section{Introduction}
\label{l:introduction}
Solutions to optimal power flow (OPF) problems provide minimum cost operating points for electric power systems in terms of a specified objective function, subject to equality constraints dictated by the non-linear power flow equations and inequalities representing engineering limits. OPF problems are generally NP-Hard~\cite{lavaei_tps,bienstock2015nphard}, even for tree networks~\cite{pascal2016nphard}, and may have multiple local optima~\cite{bukhsh_tps}. A wide variety of convex relaxation techniques have recently been applied to OPF problems~\cite{lavaei_tps,molzahn_holzer_lesieutre_demarco-large_scale_sdp_opf,dallanese2013,low_tutorial,bienstock2014,coffrin2015qc,sun2015,ghaddar2015}. Convex relaxations lower bound the optimal objective value,  can certify problem infeasibility, and, for some OPF problems, provide the global solution.

Recognizing that OPF problems are \emph{polynomial optimization} problems facilitates the application of ``moment relaxations'' from the Lasserre hierarchy~\cite{lasserre_book}. Moment relaxations generalize the semidefinite programming (SDP) relaxation proposed in~\cite{lavaei_tps} in order to solve a broader class of OPF problems at the computational cost of larger SDPs~\cite{pscc2014,cedric_msdp,ibm_opf,molzahn_hiskens-sparse_moment_opf,pscc2016,josz_molzahn-complex_hierarchy}. With increasing order in the Lasserre hierarchy, the moment relaxations provably converge to the global optima of a class of polynomial optimization problems that satisfy certain technical conditions~\cite{lasserre_book}, which includes OPF problems~\cite{josz_duality}.

Recent work~\cite{pscc2014,illustrative_example} has investigated the ability of second- and third-order moment relaxations to approach the convex hulls of several OPF problems for which the first-order relaxation (or, equivalently, the SDP relaxation of~\cite{lavaei_tps}) has a non-zero \emph{relaxation gap} (i.e., there is a gap between the global solution to the non-convex problem and the solution to the first-order relaxation). For these problems, approaching the convex hulls of the OPF problems' non-convex feasible spaces enables the relaxations to find the global optima.

Some OPF problems have objective functions that are non-linear in the power generation (e.g., convex quadratic generation costs). The unconstrained minimum of such an objective function may lie inside the convex hull of the feasible space defined by the OPF problem's constraints, but outside of the feasible space itself. For such cases, tightening the constraints of convex relaxations to more closely approximate the convex hull of the OPF problem's constraints is \emph{not sufficient} to close the relaxation gap.

In~\cite{illustrative_example}, example OPF problems using the three-bus system from~\cite{iscas2015} illustrate the ability of the moment relaxations to closely approximate the convex hulls of non-convex feasible spaces. Importantly, the moment relaxations also globally solved cases where a quadratic generation cost function was minimized at a point inside the convex hull of the feasible space, but outside the feasible space itself. In~\cite{illustrative_example}, this ability was incorrectly attributed to the non-linear relationship between the power injections and voltage phasors. (In fact, the moment relaxations have a linear relationship between the \emph{lifted variables} and the expressions for power generation.)

This paper uses a modified version of the three-bus test case from~\cite{iscas2015} to show how the constraints and objective in the moment relaxations interact. This interaction enables global solution (zero relaxation gap) of problems where the objective has an unconstrained minimizer at an infeasible point inside the convex hull of the feasible space. In contrast to constraint-tightening methods, e.g.,~\cite{coffrin2015qc,sun2015,coffrin_tightening,chen_tightening}, the second-order moment relaxation globally solves this problem.

This paper is organized as follows. Section~\ref{l:opf_overview} overviews the OPF problem. Section~\ref{l:moment} describes the Lasserre hierarchy of moment relaxations in the context of the OPF problem. Section~\ref{l:example} presents this paper's main contribution: an illustration showing how the second-order moment relaxation globally solves an OPF problem where the objective function has its minimum at an infeasible point inside the convex hull of the constraints' feasible space. Section~\ref{l:conclusion} concludes the paper.

\section{Optimal Power Flow Overview}
\label{l:opf_overview}
Consider an $n$-bus power system, where $\mathcal{N} := \left\lbrace 1, \ldots, n \right\rbrace$ is the set of all buses and $\mathcal{G}$ is the set of generator buses. Let $P_{Di} + \mathbf{j} Q_{Di}$ represent the active and reactive load demand at each bus $i \in \mathcal{N}$. Let $V_i := V_{di} + \mathbf{j} V_{qi}$ represent the voltage phasors in rectangular coordinates at each bus $i \in \mathcal{N}$. Superscripts ``max'' and ``min'' denote specified upper and lower limits. Buses without generators have maximum and minimum generation set to zero (i.e., $P_{Gi}^{\max} = P_{Gi}^{\min} = Q_{Gi}^{\max} = Q_{Gi}^{\min} = 0, \;\; \forall i\in \mathcal{N}\setminus\mathcal{G}$). Let $\mathbf{Y} := \mathbf{G} + \mathbf{j} \mathbf{B}$ denote the network admittance matrix.

The active power and reactive power generated at bus~$i$, $P_{Gi} := f_{Pi}\left(V_d,V_q\right)$ and $Q_{Gi} := f_{Qi}\left(V_d,V_q\right)$, respectively, are related to the voltages through the power flow equations:
\begin{subequations}
\small
\vspace{-5pt}
\begin{align}
\label{opf_Pbalance}
f_{Pi}\! := & P_{Di}\! +\! \sum_{k=1}^n V_{di} \left( \mathbf{G}_{ik} V_{dk}\! -\! \mathbf{B}_{ik} V_{qk} \right)\! +\! V_{qi} \left( \mathbf{B}_{ik}V_{dk}\! +\! \mathbf{G}_{ik}V_{qk} \right)  \\[-5pt]
\label{opf_Qbalance}
f_{Qi}\! := & Q_{Di}\! +\! \sum_{k=1}^n V_{qi} \left( \mathbf{G}_{ik} V_{dk}\! -\! \mathbf{B}_{ik} V_{qk}\right)\! -\! V_{di} \left( \mathbf{B}_{ik}V_{dk}\! +\! \mathbf{G}_{ik} V_{qk}\right).
\end{align}
\end{subequations}

\vspace{-10pt}

%
\noindent The squared voltage magnitude at bus~$i$ is
\begin{align}\label{Vsqr}
f_{Vi}\left(V_d,V_q\right) := V_{di}^2 + V_{qi}^2.
\end{align}
Each generator has a convex quadratic generation cost:
\begin{align}\label{cost}
f_{Ci}\left(V_d,V_q\right) = c_{2,k} \left(f_{Pi}\left(V_d,V_q\right)\right)^2 + c_{1,k} f_{Pi}\left(V_d,V_q\right) + c_{0,k}.
\end{align}

\vspace{-12pt} The OPF problem formulation considered in this paper is
%
\begin{subequations}
\label{opf}
\begin{align}
\label{opf_obj} & \min_{V_d,V_q}\quad  \sum_{i \in \mathcal{G}} f_{Ci}\left(V_d,V_q\right) \qquad \mathrm{subject\; to} \hspace{-20pt} & \\
\label{opf_P} &  \quad P_{Gi}^{\mathrm{min}} \leq f_{Pi}\left(V_d,V_q\right) \leq P_{Gi}^{\mathrm{max}} & \forall i \in \mathcal{N} \\
\label{opf_Q} &  \quad Q_{Gi}^{\mathrm{min}} \leq f_{Qi}\left(V_d,V_q\right) \leq Q_{Gi}^{\mathrm{max}} &  \forall i \in \mathcal{N} \\
\label{opf_V} &  \quad \left(V_{i}^{\mathrm{min}}\right)^2 \leq f_{Vi}\left(V_d,V_q\right) \leq \left(\vphantom{V_{i}^{\mathrm{min}}} V_{i}^{\mathrm{max}}\right)^2 &  \forall i \in \mathcal{N}  \\
\label{opf_Vref} & \quad V_{q1} = 0.
\end{align}
\end{subequations}

\vspace{-5pt}
\noindent Constraint~\eqref{opf_Vref} sets the reference bus angle to zero. 

\section{Moment Relaxations}
\label{l:moment}
The moment relaxations are developed by applying the Lasserre hierarchy of SDPs~\cite{lasserre_book} to the OPF problem~\eqref{opf}~\cite{pscc2014,cedric_msdp,ibm_opf,molzahn_hiskens-sparse_moment_opf,pscc2016}. 
%
%
Development of the moment relaxations requires several definitions. Define the vector of real decision variables $x \in \mathbb{R}^{2n}$ as $x := \begin{bmatrix} V_{d1} & \ldots & V_{dn} & V_{q1} & \ldots V_{qn} \end{bmatrix}^\intercal$. A monomial is defined using a vector $\alpha \in\mathbb{N}^{2n}$ of exponents: $x^\alpha := V_{d1}^{\alpha_1}V_{d2}^{\alpha_2}\cdots V_{qn}^{\alpha_{2n}}$. A polynomial is $h\left(x\right) := \sum_{\alpha \in \mathbb{N}^{2n}} h_{\alpha} x^{\alpha}$, where $h_{\alpha}$ is the real scalar coefficient corresponding to the monomial $x^{\alpha}$.

Define a linear functional $L_y\left\lbrace h\right\rbrace$ which replaces the monomials $x^{\alpha}$ in a polynomial $h\left(x\right)$ with scalar variables $y_{\alpha}$:
\begin{equation}
\label{Lreal}
L_y\left\lbrace h \right\rbrace := \sum_{\alpha \in \mathbb{N}^{2n}} h_{\alpha} y_{\alpha}.
\end{equation}
For a matrix $h\left(x\right)$, $L_y\left\lbrace h\right\rbrace$ is applied componentwise. 

Consider, e.g., the vector $x = \begin{bmatrix}V_{d1} & V_{d2} & V_{q1} & V_{q2} \end{bmatrix}^\intercal$ corresponding to the voltage components of a two-bus system, and the voltage-magnitude-constraint polynomial $h\left(x\right) = -\left(0.9\right)^2 + V_{d2}^2 + V_{q2}^2$. Then $L_y\left\lbrace h\right\rbrace = -\left(0.9\right)^2y_{0000} + y_{0200} + y_{0002}$. Thus, $L_y\left\lbrace h \right\rbrace$ converts a polynomial $h\left(x\right)$ to a linear function of the `lifted'' variables $y$.

For the order-$\gamma$ relaxation, define a vector $x_\gamma$ consisting of all monomials of the voltage components up to order $\gamma$:
\begin{align}
\nonumber
x_\gamma := & \left[ \begin{array}{ccccccc} 1 & V_{d1} & \ldots & V_{qn} & V_{d1}^2 & V_{d1}V_{d2} & \ldots \end{array} \right. \\* \label{x_d}
& \qquad \left.\begin{array}{cccccc} \ldots & V_{qn}^2 & V_{d1}^3 & V_{d1}^2 V_{d2} & \ldots & V_{qn}^\gamma \end{array}\right]^\intercal.
\end{align}

The moment relaxations enforce positive semidefinite constraints on so-called \emph{moment} and \emph{localizing} matrices. The symmetric moment matrix $\mathbf{M}_{\gamma}\left\lbrace y \right\rbrace$ is composed of entries $y_\alpha$ corresponding to all monomials $x^{\alpha}$ up to order $2\gamma$:
\begin{equation}
\label{eq:real_moment}
\mathbf{M}_\gamma \left\lbrace y \right\rbrace := L_y\left\lbrace x_\gamma^{\vphantom{\intercal}} x_\gamma^\intercal\right\rbrace.
\end{equation}

Symmetric localizing matrices are defined for each constraint of~\eqref{opf}. For a polynomial constraint $h\left(x\right) \geq 0$ of degree $2\eta$, the localizing matrix is:
\begin{equation}
\label{real_local}
\mathbf{M}_{\gamma - \eta} \left\lbrace h\, y \right\rbrace := L_y \left\lbrace h\, x_{\gamma-\eta}^{\vphantom{\intercal}} x_{\gamma-\eta}^{\intercal} \right\rbrace.
\end{equation}
See~\cite{pscc2014,illustrative_example} for example moment and localizing matrices.

The order-$\gamma$ moment relaxation of the OPF problem~\eqref{opf} is
\vspace{-15pt}
\begin{subequations}
\label{msosr}
\begin{align}
\label{msosr_obj}& \min_{y}\quad  L_y\left\lbrace \sum_{i \in \mathcal{G}} f_{Ci} \right\rbrace \qquad \mathrm{subject\; to} \hspace{-150pt} &  \\
\label{msosr_Pmin} &  \mathbf{M}_{\gamma-1}\left\lbrace
\left(f_{Pi} - P_i^{\min}\right) y \right\rbrace \succeq 0 &
\forall i\in\mathcal{N} \\
\label{msosr_Pmax} &  \mathbf{M}_{\gamma-1}\left\lbrace \left(P_i^{\max} - f_{Pi} \vphantom{P_i^{\min}}\right) y \right\rbrace \succeq 0 & \forall i\in\mathcal{N}\\
\label{msosr_Qmin} &  \mathbf{M}_{\gamma-1}\left\lbrace \left(f_{Qi} - Q_i^{\min}\right) y \right\rbrace \succeq 0 & \forall i\in\mathcal{N}\\
\label{msosr_Qmax} &  \mathbf{M}_{\gamma-1}\left\lbrace \left(Q_i^{\max} - f_{Qi}  \vphantom{P_i^{\min}}\right) y \right\rbrace \succeq 0 & \forall i\in\mathcal{N}\\
\label{msosr_Vmin} &  \mathbf{M}_{\gamma-1}\left\lbrace \left(f_{Vi} - \left(V_i^{\min}\right)^2\right) y \right\rbrace \succeq 0 & \forall i\in\mathcal{N}\\
\label{msosr_Vmax} &  \mathbf{M}_{\gamma-1}\left\lbrace \left(\left(V_i^{\max}\right)^2 - f_{Vi}  \vphantom{P_i^{\min}}\right) y \right\rbrace \succeq 0 & \forall k\in\mathcal{N} \\
\label{eq:msosr_Msdp} & \mathbf{M}_\gamma \left\{y\right\} \succeq 0 & \\
\label{eq:msosr_y0} & y_{00\ldots 0} = 1 & \\
\label{eq:msosr_Vref} & y_{\star\star\ldots\star\rho\star\ldots\star} = 0 & \rho = 1,\ldots,2\gamma.
\end{align}
\end{subequations}
where $\succeq 0$ indicates that the corresponding matrix is positive semidefinite and $\star$ represents any integer in $\left[ 0,\; 2\gamma -1\right]$. The constraint~\eqref{eq:msosr_y0} enforces the fact that $x^{0} = 1$. The constraint~\eqref{eq:msosr_Vref} corresponds to the angle reference $V_{q1} = 0$; the $\rho$ in \eqref{eq:msosr_Vref} is in the index $n+1$, which corresponds to $V_{q1}$.\footnote{The angle reference constraint~\eqref{opf_Vref} can alternatively be used to eliminate all terms corresponding to $V_{q1}$ to reduce the size of the matrices.}

The order-$\gamma$ relaxation yields a single global solution upon satisfaction of the rank condition

\vspace{-8pt}
\begin{equation}\label{rankcondition}
\mathrm{rank}\left(\mathbf{M}_{\gamma}\left\{y\right\}\right) = 1.
\end{equation} 
\vspace{-9pt}

\noindent The global solution $V^\ast$ to the OPF problem~\eqref{opf} is then determined by a spectral decomposition of the diagonal block of the moment matrix corresponding to the second-order monomials (i.e., $\left|\alpha\right| = 2$, where $\left|\;\cdot\; \right|$ indicates the one-norm). Specifically, let $\eta$ be a unit-length eigenvector corresponding to the non-zero eigenvalue $\lambda$ of $\left[\mathbf{M}_1 \{y\}\right]_{\left(2:k,2:k\right)}$, where $k =2n+1$ and subscripts indicate the vector entries in MATLAB notation. Then the vector $V^\ast = \sqrt{\lambda} \left(\eta_{1:n} + \mathbf{j} \eta_{\left(n+1\right):2n}\right)$ is the globally optimal voltage phasor vector.

Generally, the relaxation order $\gamma$ must be greater than or equal to half the highest degree among all objective and constraint polynomials. For OPF problems, the quadratic cost of active power generation yields a quartic polynomial in $V_d$ and $V_q$. While this suggests that $\gamma \geq 2$, the first-order relaxation (i.e., $\gamma = 1$) is formulated by minimizing $\sum_{i\in\mathcal{G}} \omega_i$, where $\omega_i$ is an auxiliary variable defined for each generator~$i\in\mathcal{G}$, and adding the second-order cone constraints 

\vspace{-18pt}
\begin{align}
\nonumber &  \left(1-c_{1,i}\,L_y\left\lbrace f_{Pi} \right\rbrace-c_{0,i} + \omega_i \right) \\ \label{msosr_quadcost} & \quad \geq \left|\left| \begin{bmatrix} \left(1+c_{1,i}\,L_y\left\lbrace f_{Pi} \right\rbrace+c_{0,i}-\omega_i \right) \\ 2\sqrt{c_{2,i}}\, L_y\left\lbrace f_{Pi} \right\rbrace \end{bmatrix} \right|\right|_2  & \forall i \in \mathcal{G}
\end{align}

\vspace{-5pt}
\noindent where $\left|\left|\,\cdot\,\right|\right|_2$ denotes the vector two-norm.

\section{Illustrative Example}
\label{l:example}
This section builds on~\cite{illustrative_example} using an OPF problem adopted from the three-bus system in~\cite{iscas2015}. The feasible space for this problem is non-convex. The chosen objective function has its unconstrained minimum at an infeasible point in the interior of the constraints' convex hull. In contrast to other relaxations, the second-order moment relaxation globally solves this problem. Note that the moment relaxations do not apply any ``artificial'' modifications to the OPF problem (e.g., the moment relaxations do not use ``penalty'' terms).

\begin{figure}[t]
\centering
\begin{tikzpicture}[scale=0.73, transform shape]

\path[draw,line width=4pt] (0,0) -- (0,2);
\draw (0,0.4) node[below right] {$1$};
\draw (-1.6,2.3) node[right] {$V_1 = 1$};
\draw (-1.6,1.8) node[right] {$\theta_1 = 0^\circ$};
\draw (-2.8,-0.5) node[right] {$P_{G1} \in\left[300,\,1200\right]~\text{MW}$};
\path[draw,line width=2pt] (-0.4,1) -- (0,1);
\draw[line width=1] (-0.8,1) circle (0.4);

\path[draw,line width=2pt] (0,1.3) -- (5,1.3);
\draw (1.6,2.4) node[below] {$R_{12} + \mathbf{j} X_{12}$};
\draw (2.5,1.8) node[below] {$= 0.15 + \mathbf{j} 0.1$};

\path[draw,line width=4pt] (5,0) -- (5,2);
\draw (5,0.4) node[below left] {$2$};
\path[draw,line width=2pt] (5,1) -- (5.4,1);
\draw[line width=1] (5.8,1) circle (0.4);
\draw (5.2,2.3) node[right] {$V_2 = 1.3$};
\draw (5.2,1.8) node[right] {$P_{G2} \in \left[0,\,50\right]~\text{MW}$};
\draw[-,line width=2pt]  (5,0.3) -- (5.3,0.3);
\draw[->,line width=2pt] (5.3,0.3) -- (5.3,-0.7);
\draw (5.4,0) node[right] {$P_{D2} + \mathbf{j} Q_{D2}$};
\draw (5.5,-0.6) node[right] {$= 30~\text{MW} + \mathbf{j} 0~\text{MVAr}$};

\path[draw,line width=2pt] (0,0.65) -- (0.5,0.65) -- (1.75,-1.5) -- (1.75,-2);
\draw (-0.2,-0.8) node[below] {$R_{13} + \mathbf{j} X_{13}$};
\draw (0.3,-1.2) node[below] {$= 0.1 + \mathbf{j} 0.05$};

\path[draw,line width=4pt] (1.5,-2) -- (3.5,-2);
\draw (1.5,-2) node[left] {$3$};
\draw[->,line width=2pt] (2.5,-2) -- (2.5,-3);
\draw[line width=1] (5.8,1) circle (0.4);
\draw (2.7,-2.3) node[right] {$P_3 = 0$};
\draw (2.7,-2.8) node[right] {$Q_3 = 0$};

\path[draw,line width=2pt] (5,0.65) -- (4.5,0.65) -- (3.25,-1.5) -- (3.25,-2);
\draw (5,-0.8) node[below] {$R_{23} + \mathbf{j} X_{23}$};
\draw (5.3,-1.2) node[below] {$= 0.001 + \mathbf{j} 0.05$};
\end{tikzpicture}
\vspace{-20pt}
\caption{Three-Bus System Adopted From~\cite{iscas2015}}
\label{f:threebussystem}

\vspace{-12pt}
\end{figure}
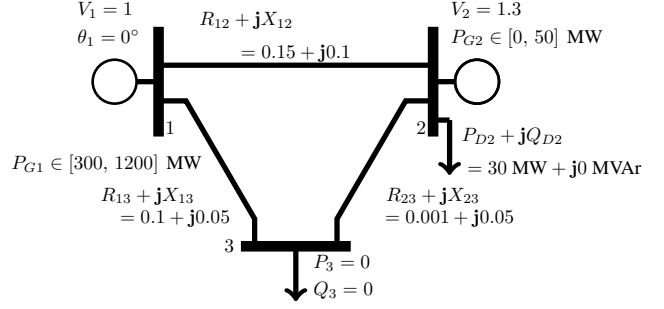

Fig.~\ref{f:threebussystem} shows the one-line diagram for the three-bus system. Voltages and line parameters are given in per~unit. The generators have no reactive power limits. Table~\ref{t:costcoefficients} gives the coefficients for the generators' quadratic cost functions. The resulting objective equals $\left(P_{G1} - 650\text{ MW}\right)^2 + 500\left(P_{G2} - 35\text{ MW}\right)^2$, which has an unconstrained minimizer at $\left(P_{G1},P_{G2}\right) = \left(650,35\right)$~MW. 

\begin{table}[t]
\centering
\caption{Generator Cost Functions}
\label{t:costcoefficients}
\begin{tabular}{|c|c|c|c|}
\hline 
Bus  & $c_{2,k}$ $(\mathrm{\$/MWh}^2)$ $\vphantom{\left(\mathrm{\$/MWh}^2\right)}$ & $c_{1,k}$ $\left(\mathrm{\$/MWh}\right)$ & $c_{0,k}$ $\left(\mathrm{\$/hr}\right)$\\ \hline\hline
1 & 1 & -1300 & 422500 \\ \hline
2 & 500 & -35000 & 612500 \\ \hline
\end{tabular}
\vspace{-8pt}
\end{table}

\begin{figure*}[t]
\centering
\subfloat[First-Order Relaxation's Feasible Space for the Three-Bus System]{\includegraphics[totalheight=0.22\textheight]{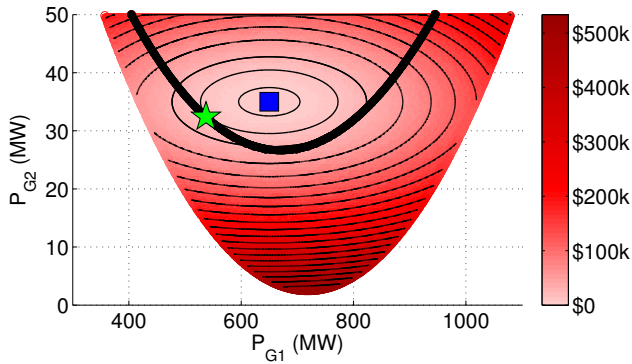}\label{f:firstorder}}
{\subfloat[Second-Order Relaxation's Feasible Space for the Three-Bus System]{\includegraphics[totalheight=0.22\textheight]{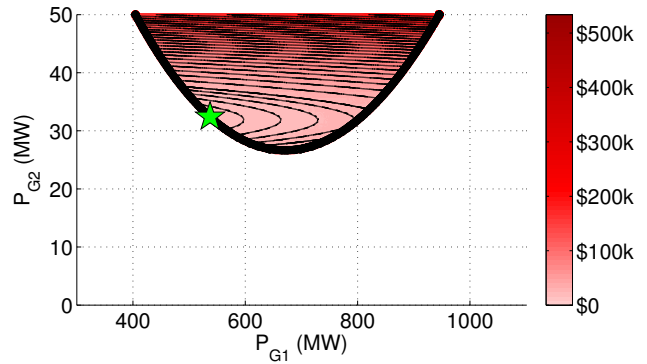}\label{f:secondorder}}} \\
\caption{A projection of the feasible spaces, in terms of active power generation, for the first- and second-order moment relaxations of the three-bus system in Fig.~\ref{f:threebussystem}. The colors represent the generation cost, with the light black lines showing cost contours. The feasible space of the OPF problem~\eqref{opf} is shown by the heavy black line. The global solution to the OPF problem~\eqref{opf} corresponds to the green star. The first-order relaxation's solution at the blue square has a non-zero relaxation gap.  Conversely, the interaction of the higher-order terms in the constraints and objective function of the second-order relaxation results in the least-cost point in the relaxation's feasible space being the global solution to the OPF problem~\eqref{opf} at the green star.}
\label{f:feasiblespace}
\vspace{-7pt}
\end{figure*}

The colored regions in Figs.~\ref{f:firstorder} and~\ref{f:secondorder} show the feasible spaces for the first- and second-order moment relaxations, with the colors representing the generation cost and the light black lines representing constant cost contours. The feasible space for the non-convex OPF problem is denoted by the heavy black curve and the global solution is at the green star.

These feasible spaces were constructed using a grid in the $P_{G1}$--$P_{G2}$ plane with a $0.5$~MW spacing. At each grid point, the first- and second-order relaxations (with objective coefficients given in Table~\ref{t:costcoefficients}) were solved with the additional constraints that the first-order relaxation's representation of the power flow equations was consistent with the specified grid point (i.e., for the grid point $\left(P_{G1},P_{G2}\right) = \left( P_{G1}^\circ,P_{G2}^\circ\right)$, both the first- and second-order moment relaxations were solved after augmenting with the constraints $\mathbf{M}_0\left\lbrace \left(f_{Pi} - P_{Gi}^\circ \right) y \right\rbrace = L_y\left\lbrace f_{Pi} - P_{Gi}^\circ \right\rbrace = 0$, $i = 1,2$). The cost shown by the colors in Figs.~\ref{f:firstorder} and~\ref{f:secondorder} corresponds to $\sum_{i\in\mathcal{G}}\omega_i$ for the first-order relaxation (cf~\eqref{msosr_quadcost}) and $\sum_{i\in\mathcal{G}}L_y\left\lbrace f_{Ci} \right\rbrace$ for the second-order relaxation. Note that the second-order relaxation's cost function~\eqref{msosr_obj} depends on terms with $\left|\alpha\right| = 2$, which appear in both the first- and second-order relaxations, and terms with $\left|\alpha\right| = 4$, which appear only in the second-order relaxation.

The solution to the first-order relaxation occurs at the blue square at $\left(P_{G1},P_{G2}\right) = \left(650,\,35\right)$~MW in Fig.~\ref{f:firstorder}, which does not match the global solution at the green star at $\left(P_{G1},P_{G2}\right) = \left(537.2,\,32.4\right)$~MW. In contrast to the first-order relaxation, Fig.~\ref{f:secondorder} shows that the second-order relaxation yields the global solution at the green star.

Observe that the unconstrained minimizer of the objective function is at the blue square in Fig.~\ref{f:firstorder}. This point is within the convex hull of the OPF problem's constraints (i.e., the convex hull of the black curve), but is not in the feasible space itself (i.e., it is not \emph{on} the black curve). A convex relaxation must enclose the feasible space of the non-convex problem. Thus, approaches for tightening the constraints of a convex relaxation can at best obtain a feasible space that matches the convex hull of the non-convex problem's constraints. If the objective function has an unconstrained minimum at an infeasible point that is in the convex hull of the non-convex problem's constraints, tightening the constraints cannot (by itself) yield a solution with zero relaxation gap. 

This has the following interpretation in the context of the three-bus problem. The solution to the first-order relaxation at the blue square in Fig.~\ref{f:firstorder} has a non-zero relaxation gap to the global optimum at the green star. Tightening the relaxation's constraints (using, e.g., the techniques in~\cite{coffrin2015qc,sun2015,coffrin_tightening,chen_tightening}) can pare down the resulting feasible space toward the convex hull of the non-convex problem's feasible space (i.e., the convex hull of the black curve in Fig.~\ref{f:firstorder}).\footnote{Applying bound tightening~\cite{coffrin_tightening} to the first-order moment relaxation augmented with the QC relaxation~\cite{coffrin2015qc} yields a relaxation whose feasible space closely matches the convex hull of this OPF problem's constraints.} However, constructing a relaxation that achieves the convex hull of the OPF problem's feasible space is not sufficient for globally solving this problem since the relaxation can still choose the point at the blue square rather than the green star in Fig.~\ref{f:firstorder}.

Focusing on the constraints alone is not sufficient for this problem; the objective function must be considered. One relevant technique adds a penalty term to the objective in an attempt to obtain a feasible solution~\cite{lavaei_mesh,lavaei_contingency}. The approach in~\cite{lavaei_mesh} penalizes the total reactive power generation.\footnote{The resulting SDP is no longer a relaxation of the original OPF problem. A solution to the penalized problem which satisfies~\eqref{rankcondition} is a \emph{feasible} point for the original problem whose worst-case optimality gap can be calculated using the lower bound from a relaxation.} Penalization approaches require the choice of penalty coefficients. It is not obvious how to obtain appropriate coefficient values for all problems, and the best known penalty coefficient ($33.76\times 10^3$~\$/(MVAr-hr)) for the three-bus problem results in a feasible point with an optimality gap of 11.8\%.\footnote{No penalization tested with the related more general approach in~\cite{lavaei_contingency} was found to substantially improve on this result.}

%
%


The second-order moment relaxation globally solves the OPF problem in Fig.~\ref{f:threebussystem} without using a penalty term or any other ``artificial'' modifications. This capability is understood in the interaction of the constraints and the objective of the second-order relaxation. Fig.~\ref{f:secondorder} shows the feasible space of the second-order relaxation, which appears to match the convex hull of the OPF problem's constraints. Importantly, despite the fact that the cost function of the OPF problem has an unconstrained minimizer in the interior of the convex hull of the constraints (i.e., the blue square in Fig.~\ref{f:firstorder}), the cost function of the second-order relaxation in Fig.~\ref{f:secondorder} is minimized at the OPF problem's global optimum at the green star.

The points inside the convex hull of the feasible space are feasible in the second-order moment relaxation. (The rank condition~\eqref{rankcondition} is not satisfied, but there exist higher-rank matrices which are feasible in~\eqref{msosr} that yield the corresponding power injections inside the convex hull.) The key observation, which can be seen in the colors and contours in Fig.~\ref{f:secondorder}, is that the higher-order terms (i.e., terms with $\left|\alpha\right| = 4$) in the objective function~\eqref{msosr_obj} result in a \emph{high cost} for the points that are infeasible in the original OPF problem~\eqref{opf}. The interaction of the second-order relaxation's constraints and objective function (i.e., the functions involving the higher-order terms with $\left|\alpha\right| = 4$) result in the global solution at the green star being the lowest cost point within the second-order relaxation's feasible space. Thus, the second-order relaxation can globally solve problems for which the objective function is minimized at an infeasible point in the convex hull of the OPF problem's constraints. Relaxations which focus solely on the constraints cannot be exact for such problems.

\vspace{-5pt}
\section{Conclusion}
\label{l:conclusion}
This paper illustrated the capabilities of moment relaxations using a small OPF problem for which the objective function has its unconstrained minimum at an infeasible point in the interior of the convex hull of the OPF problem's constraints. The second-order moment relaxation's ability to globally solve this problem is understood through the interaction of the relaxation's constraints and objective function.

\newpage
\bibliographystyle{IEEEbib}
\bibliography{refs}

\end{document}